\documentclass[12pt,reqno]{amsart}

\usepackage{color}
\usepackage{amsmath, amsthm, amssymb}
\usepackage{amsfonts}
\usepackage[ansinew]{inputenc}
\usepackage[dvips]{epsfig}
\usepackage{graphicx}
\usepackage[english]{babel}
\usepackage{enumerate}
\usepackage{hyperref}
\theoremstyle{plain}
\newtheorem{thm}{Theorem}[section]
\newtheorem{cor}[thm]{Corollary}
\newtheorem{lem}[thm]{Lemma}
\newtheorem{prop}[thm]{Proposition}

\theoremstyle{definition}

\theoremstyle{remark}

\numberwithin{equation}{section}

\newcommand{\de}{\partial}

\newcommand{\R}{\mathbb{R}}

\newcommand{\N}{\mathbb{N}}

\newcommand{\average}{{\mathchoice {\kern1ex\vcenter{\hrule height.4pt
width 6pt depth0pt} \kern-9.7pt} {\kern1ex\vcenter{\hrule
height.4pt width 4.3pt depth0pt} \kern-7pt} {} {} }}

\def\R{\mathbb{R}}

\begin{document}

\title[$C^{1,\alpha}$ estimates for the fully nonlinear Signorini problem]{$C^{1,\alpha}$ estimates for the fully nonlinear \\Signorini problem}

\author{Xavier Fernández-Real}

\address{University of Texas at Austin, Department of Mathematics, 2515 Speedway, TX 78712 Austin, USA}

\email{xevifrg@math.utexas.edu}


\thanks{The author is supported by a fellowship from ``Obra Social la Caixa''.}

\begin{abstract}
We study the regularity of solutions to the fully nonlinear thin obstacle problem. We establish local $C^{1,\alpha}$ estimates on each side of the smooth obstacle, for some small $\alpha > 0$.

Our results extend those of Milakis-Silvestre \cite{MS08} in two ways: first, we do not assume solutions nor operators to be symmetric, and second, our estimates are local, in the sense that do not rely on the boundary data.

As a consequence, we prove $C^{1,\alpha}$ regularity even when the problem is posed in general Lipschitz domains.
\end{abstract}
\maketitle
\section{Introduction}

The aim of this work is to study the regularity of the solutions to the Signorini or thin obstacle problem for fully nonlinear operators. 

Given a domain $D\subset \R^n$, the thin obstacle problem involves a function $u:D\to \R$, an obstacle $\varphi:S\to \R$ defined on a $(n-1)$-dimensional manifold $S$, a Dirichlet boundary condition given by $g:\de D\to \R$, and a second order elliptic operator $L$,
\begin{equation}
  \label{eq.classicobst}
  \left\{ \begin{array}{rcll}
  L u&=&0& \textrm{ in } D \setminus \{x\in S : u(x) = \varphi(x)\}\\
  L u&\leq&0& \textrm{ in } D\\
  u&\geq&\varphi& \textrm{ on } S \\
  u&=&g& \textrm{ on }\de D. \\
  \end{array}\right.
\end{equation}

Intuitively, one can think of it as finding the shape of a membrane with prescribed boundary conditions considering that there is a very thin obstacle forcing the membrane to be above it.

When $L$ is the Laplacian, the $C^{1,\alpha}$ regularity of solutions was first proved in 1979 by Caffarelli in \cite{Caf79}. Later, the optimal value of $\alpha$ was found by Athanasopoulos and Caffarelli in \cite{AC04}, where solutions were proved to be in $C^{1,\frac{1}{2}}$ on either side of the obstacle. More recently, this has been extended to linear operators with $x$ dependence $L = \sum a_{ij}(x)\de_{ij}u$ in \cite{Gui09, GS14, KRS15}. 

Here, we study a nonlinear version of problem \eqref{eq.classicobst}. More precisely, we study \eqref{eq.classicobst} with $Lu = F(D^2u)$, a convex fully nonlinear uniformly elliptic operator. Since all of our estimates are of local character, we consider the problem in $B_1$,
\begin{equation}
  \label{eq.thinobst}
  \left\{ \begin{array}{rcll}
  F(D^2u)&=&0& \textrm{ in } B_1\setminus \{u=\varphi\}\\
  F(D^2u)&\leq&0& \textrm{ in } B_1\\
  u&\geq&\varphi& \textrm{ on } B_1\cap\{x_n=0\}. \\
  \end{array}\right.
\end{equation}
Here, $\varphi : B_1\cap\{x_n=0\}\to \R$ is the obstacle, and we assume that it is $C^{1,1}$. We study the regularity of solutions on either side of the obstacle.

We assume that 
\begin{align}
\label{eq.F}
& F \textrm{ is convex, uniformly elliptic}\\
&\nonumber \textrm{with ellipticity constants } 0 < \lambda \leq \Lambda, \textrm{ and with } F(0)= 0.
\end{align}

When $u$ is symmetric, this problem was studied by Milakis and Silvestre in \cite{MS08}, and is equivalent to
\begin{equation}
  \label{eq.thinobstMS}
  \left\{ \begin{array}{rcll}
  F(D^2u)&=&0& \textrm{ in } B_1^+\\
  \max\{u_{x_n},\varphi-u\}&=&0& \textrm{ on } B_1\cap \{x_n  =0\}.\\
  \end{array}\right.
\end{equation}
Moreover, they also implicitly assume a symmetry condition on the operator $F$, in particular, that $F(A) = F(\tilde A)$, where $\tilde A_{in}= \tilde A_{ni} = -A_{in} = -A_{ni}$ for $i < n$ and $\tilde A_{ij} = A_{ij}$ otherwise. Under this assumption, they proved interior $C^{1,\alpha}$ regularity up to the obstacle on either side by also assuming that $ u\geq \varphi+\varepsilon$ on $\de B_1 \cap \{x_n=0\}$, for some $\varepsilon > 0$. Equivalently, they assume that the coincidence set is contained in some ball $B_{1-\delta}$ for some $\delta > 0$. This assumption is important in \cite{MS08} to prove semiconvexity of solutions. 

Our main result, Theorem~\ref{thm.main.1} below, extends the result of \cite{MS08} in two ways. First, we do not assume anything on the boundary data, so that we give a local estimate. Second, we consider also non-symmetric solutions $u$ to \eqref{eq.thinobst} with operators not necessarily satisfying any symmetry assumption, and prove $C^{1,\alpha}$ regularity for such solutions. 

In the linear case, one can symmetrise solutions to \eqref{eq.thinobst}, and then the study of such solutions reduces to problem \eqref{eq.thinobstMS}. However, in the present nonlinear setting an estimate for \eqref{eq.thinobstMS} does not imply one for \eqref{eq.thinobst}.

Our main result is the following, stating that any solution to \eqref{eq.thinobst} is $C^{1,\alpha}$ on either side of the obstacle, for some small $\alpha > 0$. 
\begin{thm}
\label{thm.main.1}
Let $F$ be a nonlinear operator satisfying \eqref{eq.F} and let $u$ be any viscosity solution to \eqref{eq.thinobst} with $\varphi\in C^{1,1}$. Then, $u\in C^{1,\alpha}(\overline{B_{1/2}^+})\cap C^{1,\alpha}(\overline{B_{1/2}^-})$ and,
\[
\|u\|_{C^{1,\alpha}(\overline{B_{1/2}^+})}+\|u\|_{C^{1,\alpha}(\overline{B_{1/2}^-})} \leq C \left(\|u\|_{L^\infty(B_1)} + \|\varphi\|_{C^{1,1}(B_1\cap \{x_{n} = 0\})} \right)
\]
for some constants $\alpha > 0$ and $C$ depending only on $n$, $\lambda$, and $\Lambda$. 
\end{thm}

Our proof of the semiconvexity of solutions is completely different from the one done in \cite{MS08} and follows by means of a Bernstein's technique. On the other hand, to prove the $C^{1,\alpha}$ regularity in the non symmetric case we follow \cite{Caf79, MS08}, but new ideas are needed. We define a symmetrised solution to the problem and follow the steps in \cite{Caf79} and \cite{MS08} using appropriate inequalities satisfied by the symmetrised solution. This yields the regularity of the symmetrised normal derivative at free boundary points. Then, we show that this implies the $C^{1, \alpha}$ regularity of the original function $u$ at free boundary points, by using the ideas from \cite{Caf89}. Finally, we show that the regularity of $u$ at free boundary points yields the regularity of the symmetrized normal derivative at all points on ${x_n=0}$, and that this yields the regularity of $u$ on either side of the obstacle.

As an immediate corollary it follows an estimate when the thin obstacle problem is posed in a bounded Lipschitz domain $D\subset \R^n$. 

\begin{cor}
\label{cor.main.1}
Let $D\subset\R^n$ be a bounded Lipschitz domain, and let $K\Subset D$. Let $F$ be a nonlinear operator satisfying \eqref{eq.F}. Let $\varphi:D\cap \{x_n = 0\}\to \R$ be a $C^{1,1}$ function, and let $u$ be the solution to
\begin{equation}
  \label{eq.thinobst_cor}
  \left\{ \begin{array}{rcll}
  F(D^2u)&=&0& \textrm{ in } D\setminus \{u = \varphi\}\\
  F(D^2u)&\leq&0& \textrm{ in } D\\
  u&\geq&\varphi& \textrm{ on } D\cap \{x_n=0\}\\
  u&=&g& \textrm{ on }\de D,\\
  \end{array}\right.
\end{equation}
for some $g\in C^0(\de D)$. Let $K^+ := K\cap \{x_n > 0\}$ and $K^- := K\cap \{x_n < 0\}$. Then, $u\in C^{1,\alpha}(\overline{K^+})\cap C^{1,\alpha}(\overline{K^-})$, with
\[
\|u\|_{C^{1,\alpha}(\overline{K^+})}+\|u\|_{C^{1,\alpha}(\overline{K^-})} \leq C \left(\|g\|_{L^\infty(\de D)} + \|\varphi\|_{C^{1,1}(D\cap \{x_n = 0\})} \right)
\]
for some constant $\alpha > 0$ depending only on $n$, $\lambda$, and $\Lambda$, and $C$ depending only on $n$, $\lambda$, $\Lambda$, $D$, and $K$.
\end{cor}

Let us introduce the notation that will be used throughout the work. We denote $x=(x',x_n) \in \R^{n}$ and 
\[
B_1^* := \{x'\in \R^{n-1} : (x', 0) \in B_1\}.\]
The obstacle $\varphi$ is defined on $B_1^*$ seen as a subset of $\R^{n}$, and problem \eqref{eq.thinobst} is written as
\begin{equation*}
  \left\{ \begin{array}{rcll}
  F(D^2u)&=&0& \textrm{ in } B_1\setminus \{(x',0): u(x',0)=\varphi(x')\}\\
  F(D^2u)&\leq&0& \textrm{ in } B_1\\
  u(x',0)&\geq&\varphi(x')& \textrm{ for } x'\in B_1^*. \\
  \end{array}\right.
\end{equation*}

We also denote 
\[
B_1^+ := \{(x', x_n)\in B_1: x_n > 0 \}, ~~~~(\de B_1)^+ = \de B_1 \cap \{x_n > 0\},
\]
and analogously we define $B_1^-$ and $(\de B_1)^-$. On the other hand, we call the coincidence set
\[
\Delta^* = \{x\in B_1^*: u(x', 0) = \varphi(x')\},~~~~\Delta = \Delta^*\times \{0\},
\] 
and its complement in $B_1^*$ is denoted by
\[
\Omega^* = B_1^*\setminus \Delta^*,~~~~\Omega = \Omega^*\times\{0\}.
\]

Our work is organised as follows. In Section~\ref{sec.2} we give a Lipschitz bound and prove semiconvexity of solutions. Then, in Section~\ref{sec.3} we prove Theorem~\ref{thm.main.1}. 

\section{Lipschitz estimate and semiconvexity}
\label{sec.2}
\subsection{Lipschitz estimate}
We begin with a proposition showing that any solution to \eqref{eq.thinobst} is Lipschitz, as long as the obstacle is $C^{1,1}$.

\begin{prop}
\label{prop.Lip}
Let $u$ be any solution to \eqref{eq.thinobst} with $F$ satisfying \eqref{eq.F} and $\varphi\in C^{1,1}$. Then $u$ is Lipschitz in $B_{1/2}$ with,
\begin{equation}
\label{eq.Lip_prop}
\|u\|_{{\rm Lip}(B_{1/2})} \leq C \left(\|u\|_{L^\infty(B_1)}+\|\varphi\|_{C^{1,1}(B_1^*)} \right),
\end{equation}
for some $C$ depending only on $n$ and the ellipticity constants $\lambda$ and $\Lambda$.
\end{prop}
\begin{proof}
We will extend the obstacle $\varphi$ to a function $h$ defined in the whole $B_1$, and we treat $u$ as a solution to a classical ``thick'' obstacle problem. We define $h$ separately in $B_1^+$ and $B_1^-$, as the solution to 
\begin{equation}
  \left\{ \begin{array}{rcll}
  F(D^2h)&=&0& \textrm{ in } B_1^+\\
  h&=&-\|u\|_{L^\infty(B_1)}& \textrm{ in } (\de B_1)^+ \\
  h(x',0)&=&\varphi(x')& \textrm{ for }x'\in B_1^*, \\
  \end{array}\right.
\end{equation}
and analogously
\begin{equation}
  \left\{ \begin{array}{rcll}
  F(D^2h)&=&0& \textrm{ in } B_1^-\\
  h&=&-\|u\|_{L^\infty(B_1)}& \textrm{ in } (\de B_1)^- \\
  h(x',0)&=&\varphi(x')& \textrm{ for }x'\in B_1^*. \\
  \end{array}\right.
\end{equation}

Notice that $h$ is Lipschitz in $B_{7/8}$; see \cite[Proposition 2.2]{MS06}. By denoting 
\[
K_0 := \|u\|_{L^\infty(B_1)} + \|\varphi\|_{C^{1,1}(B_1^*)},
\]
we have
\[
\|h\|_{{\rm Lip}(B_{7/8})} \leq CK_0,
\]
and by the maximum principle $u \geq h$. Moreover, $u$ is a solution to a classical obstacle problem in $B_1$ with $h$ as the obstacle. We show next that this implies $u$ is Lipschitz, with a quantitative estimate.

To begin with, since $h$ is Lipschitz, fixed any $x_0\in B_{1/2}$ and $0<r<1/4$, there exists some $C_0$ depending only on $n$, $\lambda$, and $\Lambda$ such that 
\begin{equation}
\label{eq.Lip1}
\sup_{B_r(x_0)} |h(x)-h(x_0)| \leq C_0 K_0 r.
\end{equation}

Notice that, by the strong maximum principle, the coincidence set $\{u = h\}$ is $\Delta$, the coincidence set of the thin obstacle problem. Suppose then that $x_0\in \Delta$, i.e., $u(x_0) = h(x_0)$. Since $u \geq h$, in particular we have that 
\begin{equation}
\label{eq.Lip2}
\inf_{B_r(x_0)} \left( u(x)-u(x_0) \right)  \geq - C_0K_0 r.
\end{equation}
because $h$ is Lipschitz. Now let
\[
q(x) = u(x)-u(x_0) +C_0K_0 r.
\]
We already know $q \geq 0$ in $B_r(x_0)$. On the other hand, from \eqref{eq.Lip1},
\[
q(x) \leq 2C_0K_0 r \textrm{~~~on~~~} B_r(x_0)\cap \Delta.
\]
Moreover, $q$ is a supersolution,
\[
F(D^2 q) = F(D^2 u) \leq 0 \textrm{~~~in~~~} B_r(x_0).
\]

Let $\bar q$ be the viscosity solution to $F(D^2 \bar q) = 0$ in $B_r(x_0)$ with $\bar q = q$ on $\de B_r(x_0)$. We have $\bar q \leq q$ in $B_r(x_0)$ and by the non-negativity of $\bar q$ on the boundary, $\bar q \geq 0$ in $B_r(x_0)$.

Thus, $q < \bar q + 2C_0K_0r$ on $\de B_r(x_0)$, and $q \leq \bar q + 2C_0K_0r$ in  $B_r(x_0)\cap \Delta$. Therefore,
\[
q \leq \bar q + 2C_0K_0r \textrm{~~~in~~~} B_r(x_0).
\]

On the other hand, we know $0\leq \bar q(x_0) \leq q (x_0) = C_0K_0r$, and by the Harnack inequality, $\bar q \leq CC_0K_0r$ in $B_{r/2}(x_0)$. Putting all together we obtain that $u(x) -u(x_0) \leq CC_0K_0r$ for some constant $C >0$. Thus, combining this with \eqref{eq.Lip2},
\begin{equation}
\label{eq.Lip3}
\sup_{B_r(x_0)} |u(x)-u(x_0)| \leq C\left(\|u\|_{L^\infty(B_1)}+\|\varphi\|_{C^{1,1}(B_{1}^*)}\right)r,
\end{equation}
for some constant $C$ depending only on $n$, $\lambda$, and $\Lambda$.

We have obtained that the solution is Lipschitz on points of the coincidence set. Let us use interior estimates to deduce Lipschitz regularity inside $B_{1/2}$. 

Take any points $x, y\in B_{1/2}$, and let $r = |x-y|$. Define
\[
\rho := \min \{\textrm{dist}(x, \Delta), \textrm{dist}(y, \Delta)\},
\]
and let $x^*, y^*\in\Delta$, $x^* = (x',0)$, $y^* = (y',0)$ for $x', y'\in \Delta^*$, be such that $\textrm{dist}(x, \Delta)= |x-x^*|$ and $\textrm{dist}(y, \Delta)= |y-y^*|$. We now separate two cases:

$\bullet$ If $\rho \leq 4r$, then
\begin{align*}
|u(x)-u(y)| & \leq |u(x)-u(x^*)|+|u(y)-u(y^*)|+|\varphi(x')-\varphi(y')|\\
& \leq C\rho + C(r+\rho)+2C(r+\rho) \leq Cr
\end{align*}
for some constant $C$. We are using here that $\varphi$ is Lipschitz and that if $|x-x^*| = \rho$, then $|y-y^*| \leq r+ \rho$ and $|x^*-y^*| \leq 2(r+\rho)$.

$\bullet$ If $\rho > 4r$, we can use interior estimates. Suppose $x$ is such that $\textrm{dist}(x,\Delta) = \rho$, and notice $B_{\rho/2}(x) \subset B_1\setminus \Delta$, so that in $B_{\rho/2}(x)$, $F(D^2 u) = 0$. We can now use the interior Lipschitz estimates (see, for example, \cite[Chapter 5]{CC95}),
\[
[u]_{{\rm Lip}(B_{\rho/4})} \leq \frac{C}{\rho}\textrm{osc}_{B_{\rho/2}(x)} u \leq C
\]
for some constant $C$. We are using here that the supremum and the infimum of $u$ in $B_{\rho/2}(x)$ are controlled respectively by $C\rho + \varphi(x^*)$ and $-C\rho + \varphi(x^*)$.

Thus, we have proved that the solution is Lipschitz in $B_{1/2}$, with the estimate \eqref{eq.Lip_prop}.
\end{proof}

\subsection{Preliminary consideration}
\label{ssec.2.2}
Before continuing to prove the semiconvexity and semiconcavity result, we introduce a change of variables that will be useful in this section and the next one. Notice that, given a function $w$, we can express the nonlinear operator $F$ as
\[
F(D^2 w(x)) = \sup_{\gamma\in \Gamma} \left( L_\gamma^{ij}\de_{x_ix_j}w(x)+c_{\gamma}
 \right),\]
for some family of symmetric uniformly elliptic operators with ellipticity constants $\lambda$ and $\Lambda$, $L_\gamma^{ij}\de_{x_ix_j}$, indexed by $\gamma\in \Gamma$. Since $F(0) = 0$, there is some symmetric uniformly elliptic operator from this family given by a matrix $\hat{L}$ such that 
\[
{\rm tr}(\hat L D^2 w(x))=\hat{L}^{ij}\de_{x_ix_j}w(x) \leq F(D^2w(x)).
\]
We now change coordinates in such a way that the matrix of this operator in the new coordinates, denoted $\hat{L}_A$, fulfils $\hat{L}_A^{in} = \hat{L}_A^{ni} = 0$ for $i < n$. More precisely, if we denote $\hat L '$ the matrix in ${\rm Sym}_{n-1}$ given by the $n-1$ first indices of $\hat L$, and we denote $\hat L_n' = (\hat L^{in})_{1\leq i \leq n-1}$ the vector of $\R^{n-1}$, we change variables as
\[
x\mapsto y = Ax,
\]
where $A$ is the matrix given by 
\[
A:= \left(
\begin{array}{c|c}
  \raisebox{-10pt}{{\large\mbox{{Id}}}}_{n-1} & \\[-2.5ex]
     & - \bar a\\ 
        & \\ \hline
  0\dots 0 & 1 \\[-0.5ex]
\end{array}
\right),
\]
and $ \bar a = (\hat L')^{-1}\cdot \hat L_n'$ is a vector in $\R^{n-1}$. We define the new nonlinear operator $\tilde{F}$ as
\[
\tilde{F}(N) = F(A^{T} N A), \textrm{ for all } N \in {\rm Sym}_n, 
\]
so that it is consistent with the change of variables, in the sense that if $\tilde{w}(y) = w(A^{-1} y)$, then $F(D^2 w(x)) = \tilde{F}(D^2\tilde{w}(y))$.

We trivially have that $\tilde{F}$ is convex and $\tilde{F} (0) = 0$. In the new coordinates we still have that $\hat L_A^{ij} \de_{y_iy_j}$ is a symmetric uniformly elliptic operator, but now the ellipticity constants $\lambda$ and $\Lambda$ have changed depending only on $n$, $\lambda$, and $\Lambda$. The same occurs with all the operators in the family defining $F$, so that after changing coordinates, $F$ is still a convex uniformly elliptic operator with ellipticity constants depending only on $n$, $\lambda$, and $\Lambda$. Indeed, for any matrices $N, N_P\in {\rm Sym}_{n-1}$ with $N_P \geq 0$ we have that (using the definition of uniform ellipticity in \cite[Chapter 2]{CC95} and noticing that $A^T N_P A \geq 0$), 
\[
 \|A^{-1}\|^{-2} \|N_P\| \leq \lambda\|A^T N_P A\| \leq \tilde{F}(N+N_P) - \tilde{F}(N) \leq \Lambda \|A^T N_P A \|\leq \Lambda \|A\|^2 \|N_P\|, 
\]
and it is easy to bound $\|A^{-1}\|$ and $\|A\|$ from the definition of $A$, depending only on $n$, $\lambda$ and $\Lambda$.

After changing variables, the regularity of the solution remains the same up to multiplicative constants in the bounds depending only on $n$, $\lambda$, and $\Lambda$.

As an abuse of notation we will call the new variables $(x',x_n)$, the new operator $F$, and the new ellipticity constants $\lambda$ and $\Lambda$, understanding that they might depend on the original ellipticity constants and the dimension, $n$. This will not be a problem, since in all the statements of the present work $n$, $\lambda$, and $\Lambda$ appear together in the dependence of the constants. 

Thus, throughout the paper we will assume that there exists a fixed symmetric uniformly elliptic operator $\hat L$ such that
\begin{equation}
\label{eq.Lhat}
\hat{L}^{ij}\de_{x_ix_j}w \leq F(D^2w),~~{\rm and}~~ \hat{L}^{in} = \hat{L}^{ni} = 0 ~~{\rm for}~~i <n.
\end{equation}

This change of variables is useful because, for any function $w$,
\[
\hat{L}^{ij}\de_{x_ix_j}(w(x',-x_n))=\hat{L}^{ij}(\de_{x_ix_j}w)(x',-x_n),
\]
which will allow us to symmetrise the solution and still have a supersolution for the Pucci extremal operator $\mathcal{M}^-$. We also use it to prove a semiconcavity result from semiconvexity in the following proof of Proposition \ref{prop.semicon}.

\subsection{Semiconvexity and semiconcavity estimates}
We next prove the semiconvexity of solutions in the directions parallel to the domain of the obstacle. To do it, we use a Bernstein's technique in the spirit of \cite{AC04}.
\begin{prop}
\label{prop.semicon}
Let $u$ be the solution to \eqref{eq.thinobst}. Then
\begin{enumerate}[(a)]
\item (Semiconvexity) If $\tau = (\tau^*,0)$, with $\tau^*$ a unit vector in $\R^{n-1}$,
\[
\inf_{B_{3/4}} u_{\tau \tau} \geq -C\left(\|u\|_{L^\infty(B_1)}+\|\varphi\|_{C^{1,1}(B_1^*)}\right),
\]
for some constant $C$ depending only on $n$, $\lambda$, and $\Lambda$.
\item (Semiconcavity) Similarly, in the direction normal to $B_1^*\times\{0\}$,
\[
\sup_{B_{3/4}} u_{x_nx_n} \leq C\left(\|u\|_{L^\infty(B_1)}+\|\varphi\|_{C^{1,1}(B_1^*)}\right),
\]
for some constant $C$ depending only on $n$, $\lambda$, and $\Lambda$.
\end{enumerate}
\end{prop}

\begin{proof}
The second part, (b), follows from (a) using the definition of uniformly elliptic operator and the fact that we changed variables (in the previous subsection) in order to have matrix $\hat L$ fulfilling \eqref{eq.Lhat}. We denote by $\hat L'$ and $D^2_{n-1}u$ the square matrices corresponding to the $n-1$ first indices of $\hat L$ and $D^2 u$ respectively. Now, from
\[
\hat{L}^{ij} \de_{x_ix_j} u(x) \leq 0,~~~\hat{L}^{in} = \hat{L}^{ni} = 0 ~\textrm{ for }~ i < n,
\]
and
\[D^2_{n-1} u \geq -C\left(\|u\|_{L^\infty(B_1)}+\|\varphi\|_{C^{1,1}(B_1^*)}\right){\rm Id}_{n-1},
\]
we directly obtain that 
\[
\hat L^{nn}\de_{x_nx_n} u \leq - \sum_{i, j = 1}^{n-1} \hat L^{ij}\de_{x_ix_j} u \leq C \left(\|u\|_{L^\infty(B_1)}+\|\varphi\|_{C^{1,1}(B_1^*)}\right){\rm tr} \hat L'.
\]
The desired bound follows because $\hat{L}^{nn}$ is bounded below by $\lambda$ and ${\rm tr}(\hat L' )$ is bounded above by $(n-1)\Lambda$. 

Let us prove (a). As in the proof of Proposition \ref{prop.Lip}, we define $h$ as the solution to
\begin{equation}
  \left\{ \begin{array}{rcll}
  F(D^2h)&=&0& \textrm{in } B_1^+\\
  h&=&-\|u\|_{L^\infty(B_1)}& \textrm{in } (\de B_1)^+ \\
  h(x',0)&=&\varphi(x')& x'\in B_1^* \\
  \end{array}\right.~~~~~~~~
  \left\{ \begin{array}{rcll}
  F(D^2h)&=&0& \textrm{in } B_1^-\\
  h&=&-\|u\|_{L^\infty(B_1)}& \textrm{in } (\de B_1)^- \\
  h(x',0)&=&\varphi(x')& x'\in B_1^*. \\
  \end{array}\right.
\end{equation} 
Recall that $h$ is Lipschitz and that, by the strong maximum principle, $u > h$ in $B_{1/2}^+$ and $B_{1/2}^-$. 

Define now, for $\varepsilon > 0$, 
\[
\bar h_\varepsilon(x', x_n) := \varphi(x')-\frac{x_n^2}{\varepsilon}\]
and
\[
h_\varepsilon(x', x_n) := \max\left\{h(x', x_n), \bar h_\varepsilon(x', x_n)\right\}.
\]
Since, $h$ is Lipschitz continuous and $h(x', 0) = \bar h_\varepsilon (x', 0)$, this implies that there exists a constant $C > 0$ depending only on $n$, $\lambda$, and $\Lambda$ such that
\begin{equation}
\label{eq.epsliph}
h(x', x_n) > \bar h_\varepsilon (x', x_n) ~~ \textrm{ for } ~~ |x_n| > CK_0\varepsilon,
\end{equation}
where we define
\[
K_0:= \|u\|_{L^\infty(B_1)}+\|\varphi\|_{C^{1,1}(B_1^*)}.
\]
In particular, $h_\varepsilon$ is Lipschitz continuous in $B_{7/8}$, uniformly on $\varepsilon$. 

Let $u_\varepsilon$ be the solution to the ``thick'' obstacle problem with obstacle $h_\varepsilon$, 
\begin{equation}
\label{eq.thickeps}
  \left\{ \begin{array}{rcll}
  F(D^2u_\varepsilon)&=&0& \textrm{ in } B_1\setminus \{u_\varepsilon = h_\varepsilon\}\\
  F(D^2u_\varepsilon)&\leq&0& \textrm{ in } B_1\\
  u_\varepsilon&=&\max\{u,\bar h_\varepsilon\}& \textrm{ on } \de (B_1^+)\\
  u_\varepsilon&\geq&h_\varepsilon& \textrm{ in } B_1^+,\\
  \end{array}\right.
\end{equation}
and the analogous expression in $B_1^-$. By \eqref{eq.epsliph}, the coincidence set satisfies 
\[
\{u_\varepsilon = h_\varepsilon\}\subset \{\bar h_\varepsilon > h\} \subset \{(x', x_n)\in B_1 : |x_n| \leq CK_0\varepsilon\}
\]
for some $C >0$. We want to bound $\de_{\tau\tau}u_{\varepsilon}$ from below independently of $\varepsilon$.

Notice that $D^2(u_\varepsilon - h_\varepsilon) \geq 0$ in the coincidence set, and since $u_\varepsilon \geq h_\varepsilon$, this also occurs along the free boundary. By the definition of $\bar h_\varepsilon$ and recalling that $h_\varepsilon = \bar h_\varepsilon$ in the coincidence set, this implies $\de_{\tau\tau}u_{\varepsilon} \geq -CK_0$ in $\{u_\varepsilon = h_\varepsilon\}\cap B_{7/8}$, for some constant $C$ depending only on $n$, $\lambda$, and $\Lambda$. Thus, it is enough to check that $\de_{\tau\tau}u_{\varepsilon}$ is uniformly bounded from below outside the coincidence set. We proceed by means of a Bernstein's technique.

Let $\eta \in C^\infty_c({B_{7/8}})$ be a smooth, cutoff function, with $0\leq \eta\leq 1$ and $\eta \equiv 1$ in $B_{3/4}$. Define 
\[
f_\varepsilon (x) = \eta(x) \de_{\tau\tau}u_{\varepsilon}(x) - \mu |\nabla u_\varepsilon(x)|^2
\]
for some constant $\mu$ to be determined later. Notice that, since $h_\varepsilon$ is Lipschitz con\-ti\-nuous independently of $\varepsilon$ in $B_{7/8}$, then $|\nabla u_\varepsilon(x)|$ is bounded independently of $\varepsilon$ in $B_{7/8}$. If the minimum $x_0$ in $B_{7/8}$ is attained in the coincidence set, then $\de_{\tau\tau}u_{\varepsilon}(x_0) \geq -CK_0$ and we get that for every $x\in B_{3/4}$, 
\begin{equation}
\label{eq.bernst1}
 \de_{\tau\tau}u_{\varepsilon}(x)\geq -CK_0 - \mu |\nabla u_\varepsilon(x_0)|^2+ \mu |\nabla u_\varepsilon(x)|^2 \geq -CK_0-\mu \|\nabla u_\varepsilon\|_{L^\infty(B_{7/8})}^2.
\end{equation}
If the minimum $x_0$ is attained at the boundary, $\de B_{7/8}$, then for every $x\in B_{3/4}$, 
\begin{equation}
\label{eq.bernst2}
 \de_{\tau\tau}u_{\varepsilon}(x)\geq -\mu |\nabla u_\varepsilon(x_0)|^2+ \mu |\nabla u_\varepsilon(x)|^2 \geq -\mu \|\nabla u_\varepsilon\|_{L^\infty(B_{7/8})}^2.
\end{equation}

Let us assume now that the minimum $x_0$ of $f_\varepsilon$ in $B_{7/8}$ is attained at some interior point $x_0$ outside the coincidence set $\{u_\varepsilon = h_\varepsilon\}$.

Let us also assume that the operator $F$ not only is convex, but also $F\in C^\infty$, so that solutions are $C^4$ outside the coincidence set (see the end of the proof for the general case $F$ Lipschitz). In this case, the linearised operator of $F$ at $x_0$,
\[
L_0 v = a_{ij} v_{ij} := F_{ij}(D^2 u_\varepsilon(x_0)) v_{ij},
\]
is uniformly elliptic with ellipticity constants $\lambda$ and $\Lambda$. Moreover, for any $\rho\in S^{n-1}$,
\begin{equation}
\label{eq.L0}
L_0 u_\varepsilon(x_0) \geq 0,~~~L_0 \de_\rho u_{\varepsilon}(x_0) = 0,~~~ L_0 \de_{\rho\rho} u_{\varepsilon}(x_0) \leq  0.
\end{equation}
This is a standard result, which can be found in \cite[Lemma 9.2]{CC95}. 

For simplicity in the following computations we denote $w = u_\varepsilon$. If $x_0$ is an interior minimum of $f_\varepsilon$ (which is a $C^2$ function) in $B_{7/8}$, then
\begin{equation}
\label{eq.ellip1}
0 = \nabla f_\varepsilon(x_0) = (\nabla\eta w_{\tau \tau}+\eta\nabla w_{\tau \tau}-2\mu w_i\nabla w_i) (x_0),
\end{equation}
and by \eqref{eq.L0} and the fact that $(a_{ij})$ is elliptic,
\begin{equation}
\label{eq.ellip2}
0 \leq a_{ij}f_{\varepsilon, ij}(x_0) \leq \left(a_{ij}\eta_{ij} w_{\tau \tau} + 2 a_{ij}\eta_i w_{\tau \tau,j} - 2 \mu a_{ij}w_{kj} w_{ki}\right)(x_0).
\end{equation}

Combining \eqref{eq.ellip1} and \eqref{eq.ellip2}, we find
\begin{equation}
\label{eq.ellip2.5}
0\leq \left(\left(a_{ij}\eta_{ij} -2\frac{a_{ij}\eta_i\eta_j}{\eta}\right)w_{\tau\tau}- 2 \mu a_{ij}w_{kj} w_{ki} + 4\frac{\mu a_{ij}\eta_i w_{kj} w_k}{\eta}\right)(x_0).
\end{equation}
Observe that $|\nabla \eta |^2 \leq C\eta$ (since $\sqrt{\eta}$ is Lipschitz). Therefore, for some constants $C_0$ and $C_1$ depending only on $n$ and $\Lambda$,
\[
0\leq \left(C_0 |w_{\tau\tau}| + \mu C_1 |D^2 w||\nabla w| - 2 \mu a_{ij}w_{kj} w_{ki}\right)(x_0).
\]
Using $|w_{\tau \tau}(x_0)| \leq |D^2w(x_0)|$ and the uniform ellipticity of $(a_{ij})$,
\[
a_{ij}w_{ki}w_{kj} \geq \lambda C(n) |D^2 w|^2,
\]
we obtain
\[
|D^2 w(x_0)| \leq  \frac{C_0}{\mu} + C_1|\nabla w(x_0)|,
\]
for some constants $C_0$ and $C_1$ depending now also on $\lambda$. Now, since $x_0$ is a minimum in $B_{7/8}$, for any $x\in B_{3/4}$,
\begin{align}
\nonumber w_{\tau\tau}(x)& \geq \eta(x_0)w_{\tau\tau}(x_0) - \mu |\nabla u_\varepsilon(x_0)|^2+ \mu |\nabla u_\varepsilon(x)|^2 \\
 \label{eq.bernst3} & \geq -|D^2w(x_0)|-\mu \|\nabla u_\varepsilon\|_{L^\infty(B_{7/8})}^2\\
\nonumber & \geq -\frac{C_0}{\mu} - C_1\|\nabla u_\varepsilon\|_{L^\infty(B_{7/8})}-\mu \|\nabla u_\varepsilon\|_{L^\infty(B_{7/8})}^2.
\end{align} 

We now fix $\mu = \|\nabla u_\varepsilon\|_{L^\infty(B_{7/8})}^{-1}$. Notice that, in all three cases \eqref{eq.bernst1}, \eqref{eq.bernst2}, and \eqref{eq.bernst3}, we reach that for some constant $C$ depending only on $n$, $\lambda$, and $\Lambda$, 
\[
\inf_{B_{3/4}} \de_{\tau\tau} u_{\varepsilon}  \geq -C\left( \sup_{B_{7/8}} |\nabla u_\varepsilon| +K_0\right).
\]
We had already seen that $u_\varepsilon$ is Lipschitz continuous independently of $\varepsilon > 0$ and controlled by the Lipschitz norm of $u$, so that by Proposition \eqref{prop.Lip},
\begin{equation}
\label{eq.ellip3}
\inf_{B_{3/4}} \de_{\tau\tau} u_{\varepsilon}  \geq -C\left( \|u\|_{{\rm Lip}(B_{7/8})} + \|\varphi\|_{C^{1,1}(B_1^*)}+K_0\right) \geq -C\left( \|u\|_{L^\infty(B_{1})} + \|\varphi\|_{C^{1,1}(B_1^*)}\right) .
\end{equation}

If $F$ is not smooth, then it can be regularised convoluting with a mollifier in the space of symmetric matrices, so that it can be approximated uniformly in compact sets by a sequence $\{F_k\}_{k\in \N}$ of convex smooth uniformly elliptic operators with ellipticity constants $\lambda$ and $\Lambda$; also, by subtracting $F_k(0)$, we can assume $F_k(0) = 0$. Note that, in $B_{7/8}$ and for every $\varepsilon > 0$ we have uniform $C^{1,\gamma}$ estimates in $k$ for the solutions to \eqref{eq.thickeps} with operators $F_k$, since the obstacle $h$ is in $C^{1,1}$ in a neighbourhood of the free boundary. By Arzelà-Ascoli there exists a subsequence converging uniformly, and therefore, the estimate \eqref{eq.ellip3} can be extended to solutions of \eqref{eq.thickeps} with operators not necessarily smooth. Thus, \eqref{eq.ellip3} follows for any $F$ not necessarily $C^\infty$.

Note that $u_\varepsilon$ converges uniformly to $u$, since for all $\delta > 0$, there exists some $\varepsilon > 0$ small enough such that $u+\delta > u_\varepsilon\geq u$ in $B_1$. 

Since the right-hand side of \eqref{eq.ellip3} is independent of $\varepsilon$, and $u_\varepsilon$ converges uniformly to $u$ in $B_{7/8}$ as $\varepsilon \downarrow 0$, we finally obtain
\begin{equation}
\label{eq.ellip4}
\inf_{B_{3/4}} u_{\tau\tau}  \geq -C\left( \|u\|_{L^\infty(B_{1})} + \|\varphi\|_{C^{1,1}(B_1^*)}\right),
\end{equation}
as desired. 
\end{proof}
\section{$C^{1,\alpha}$ estimate}
\label{sec.3}
\subsection{A symmetrised solution} By the results in the previous section we know that $\nabla u$ is bounded in the interior of $B_1$. Moreover, $u_{x_nx_n}$ is bounded from above inside $B_1$. In particular, the following limit exists
\begin{equation}
\label{eq.sigmadef}
\sigma(x') = \lim_{{x_n} \downarrow 0^+} u_{x_n}(x', x_n) - \lim_{x_n\uparrow 0^-} u_{x_n}(x',x_n) = \lim_{x_n \downarrow 0^+} \bigl( u_{x_n}(x',x_n) - u_{x_n}(x',-x_n) \bigr).
\end{equation}

A main step towards Theorem~\ref{thm.main.1} consists of proving that $\sigma\in C^\alpha(B_{1/2}^*)$ for some $\alpha  >0$. We will prove this in this section. 

We begin by noticing that $\sigma(x') = 0$ for $x'\in \Omega^*$ (by the $C^{2,\alpha}$ interior estimates), where we recall that $\Omega^* :=\{x'\in B_1^*: u(x',0)> \varphi(x')\}$. In general, however, we have the following:

\begin{lem}
\label{lem.sigmaneg}
The function $\sigma$ defined by \eqref{eq.sigmadef} is non-positive, i.e., $\sigma\leq 0 $ in $B_1^*$.
\end{lem}

\begin{proof}
Suppose it is not true, and there exists some $\bar x'\in B_1^*$ such that $\sigma(\bar x') > 0$. Let $\delta > 0$ be such that $B_\delta^*(\bar x')\subset B_1^*$, so that by the semiconcavity in Proposition~\ref{prop.semicon} applied to $B_{\delta/2}((\bar x', 0))$, $u_{x_nx_n}(\bar x', 0) \leq C$ for some constant $C$, that now depends also on $\delta$. However, 
\[
\sigma(\bar x') = \lim_{x_n \downarrow 0^+} \left( u_{x_n}(\bar x',x_n) - u_{x_n}(\bar x',-x_n) \right) > 0,
\]
which means 
\[
 \frac{u_{x_n}(\bar x',x_n) - u_{x_n}(\bar x',-x_n) }{2x_n} \to +\infty, {~~\rm as~~}x_n \downarrow 0^+,
\]
a contradiction with the bound in $u_{x_nx_n}$.
\end{proof}

We will now adapt the ideas of \cite{Caf79} to our non-symmetric setting. For this, we use a symmetrised solution, defined as follows 
\begin{equation}
\label{eq.usim}
v(x', x_n) := \frac{u(x', x_n)+u(x', -x_n)}{2},~\textrm{ for } (x',x_n)\in \overline{B_1}.
\end{equation}
Here $u$ is any solution to \eqref{eq.thinobst}.

Notice that 
\begin{equation}
\label{eq.sigmaus}
\sigma(x') = 2\lim_{x_n \downarrow 0^+} v_{x_n}(x',x_n) \leq 0
\end{equation}
is well defined, and in particular, we have that
\begin{equation}
\label{eq.sigmaus2}
\sigma(x') = 2v_{x_n}(x', 0) = 0, \textrm{ for } x'\in \Omega^*.
\end{equation}

The following result follows from the results in the previous section. We will use the notation $\mathcal{M}^+$ and $\mathcal{M}^-$ to refer to the Pucci's extremal operators with the implicit ellipticity constants $\lambda$ and $\Lambda$ (see \cite[Chapter 2]{CC95} for the definition and basic properties of such operators).

\begin{lem}
\label{lem.usim}
Let $u$ be a solution to the nonlinear thin obstacle problem \eqref{eq.thinobst}, and let $v$ be defined by \eqref{eq.usim}. Then $v$ is Lipschitz in $\overline{B_{1/2}^+}$ and satisfies
\begin{equation}
\label{eq.usproblem}
  \left\{ \begin{array}{rcll}
  \mathcal{M}^-(D^2v)&\leq&0& \textrm{ in } B_1,\\
  \max\{v_{x_n}(x',0), \varphi(x')-v(x',0)\}&=&0& \textrm{ for }x'\in B_1^*.
  \end{array}\right.
\end{equation}
Moreover, 
\begin{enumerate}[(a)]
\item (Semiconvexity) If $\tau = (\tau^*,0)$, with $\tau^*$ a unit vector in $\R^{n-1}$,
\[
\inf_{B_{3/4}} v_{\tau \tau} \geq -C\left(\|u\|_{L^\infty(B_1)}+\|\varphi\|_{C^{1,1}(B_1^*)}\right),
\]
for some constant $C$ depending only on $n$, $\lambda$, and $\Lambda$.
\item (Semiconcavity) In the direction normal to $B_1^*\times\{0\}$,
\[
\sup_{B_{3/4}} v_{x_nx_n} \leq C\left(\|u\|_{L^\infty(B_1)}+\|\varphi\|_{C^{1,1}(B_1^*)}\right),
\]
for some constant $C$ depending only on $n$, $\lambda$, and $\Lambda$.
\end{enumerate}
\end{lem}
\begin{proof}
The Lipschitz regularity comes from the Lipschitz regularity in $u$, proved in Proposition~\ref{prop.Lip}. 

In \eqref{eq.usproblem} the first inequality follows thanks to the change of variables introduced in Subsection \ref{ssec.2.2}. Indeed, there exists some operator given by a matrix $\hat L$ as in \eqref{eq.Lhat} uniformly elliptic with ellipticity constants $\lambda$ and $\Lambda$ such that 
\[
\hat L^{ij} \de_{x_ix_j} (u(x',-x_n)) = \hat L^{ij} (\de_{x_ix_j} u)(x',-x_n) \leq F((D^2 u)(x',-x_n))\leq 0,
\]
so that 
\[
\mathcal{M}^-(D^2 v ) \leq \hat L^{ij} \de_{x_ix_j} v \leq 0,
\]
as we wanted.

The second expression in \eqref{eq.usproblem} follows from equations \eqref{eq.sigmaus}-\eqref{eq.sigmaus2}, Lemma~\ref{lem.sigmaneg} and the fact that $v(x',0) = u(x',0)$ for $x'\in B_1^*$. 

Finally, the semiconvexity and semiconcavity follow from Proposition~\ref{prop.semicon}.
\end{proof}

\subsection{Regularity for $\sigma$ on free boundary points}
The next steps are very similar to those in \cite{Caf79} (and \cite{MS08}), but we adapt them to the symmetrised solution $v$ instead of $u$. For completeness, we provide all the details. We begin with the following lemma, corresponding to \cite[Lemma 2]{Caf79} (or \cite[Lemma 3.3]{MS08}).

In the next result, we call $\varphi$ the extension of the obstacle to $B_1$, i.e. $\varphi(x',x_n) := \varphi(x')$. 

\begin{lem}
\label{lem.hx0}
Let $v$ be the symmetrised solution \eqref{eq.usim}. Let $\kappa$ be a constant such that $\kappa > \sup|\varphi_{\tau\tau}|$ for any $\tau$ a unit vector in $\R^{n-1}\times\{0\}$. Let $x_0\in \Omega$ fixed and $\psi_{x_0}$ denote the function
\[
\psi_{x_0} = \varphi(x_0)+\nabla \varphi (x_0) \cdot (x-x_0) + \kappa |x-x_0|^2 - \kappa (n-1)\frac{\Lambda}{\lambda} x_n^2.
\]

Then, for any open set $U_{x_0}$ such that $x_0\in U_{x_0} \subset B_1$,
\[
\sup_{\de U_{x_0}\cap \{x_n > 0\} } (v-\psi_{x_0}) \geq 0.
\]
\end{lem}
\begin{proof}
Define $w = v-\psi_{x_0}$ and notice that by definition of $\psi_{x_0}$ and the fact that $v$ is a supersolution for $\mathcal{M}^-$, we have $w(x_0) \geq 0$ and $\mathcal{M}^-(D^2w) \leq 0$. Therefore, we can apply the maximum principle on $U_{x_0}\setminus \Delta$ (recall $\Delta$ is the coincidence set) and use the symmetry of $w$ to obtain that 
\[
\sup_{\de (U_{x_0}\setminus\Delta)\cap \{x_n \geq 0\} } (v-\psi_{x_0}) \geq 0.
\]

Now notice that on the set $\{v = u = \varphi\}$ we have that $\psi_{x_0} > \varphi$, since $x_0 \in \Omega$ and $\kappa > \sup |\varphi_{\tau \tau}|$. Thus, $v-\psi_{x_0} < 0$ on this set, so that
\[
\sup_{\de (U_{x_0}\setminus\Delta)\cap \{x_n \geq 0\} } (v-\psi_{x_0}) = \sup_{\de U_{x_0}\cap \{x_n > 0\} } (v-\psi_{x_0}) \geq 0,
\]
and we are done.
\end{proof}

We now proceed with the following lemma, corresponding to \cite[Lemma 2]{Caf79} (or \cite[Lemma 3.4]{MS08}). 

\begin{lem}
\label{lem.sgamma}
Let $v$ be the symmetrised solution as defined in \eqref{eq.usim}, and let $\sigma$ as defined in \eqref{eq.sigmadef}-\eqref{eq.sigmaus}. Let $x_0=(x_0',0)\in \Omega$ and define $S_\gamma = \{x': \sigma(x') > -\gamma\}$. Then, for suitable positive constants $C$, $\overline{C}$, and $\gamma_0$ and for all $\gamma \in (0,\gamma_0)$ there exists a ball $B_{C\gamma}^*(\bar{x}')$ for $\bar{x}'\in B_1^*$ such that 
\[
B_{C\gamma}^*(\bar{x}') \subset B_{\overline C \gamma}^* (x_0')\cap S_\gamma.
\]
The constants $C$, $\bar C$, and $\gamma_0$ depend only on $n$, $\lambda$, $\Lambda$, $\|\varphi\|_{C^{1,1}(B_1^*)}$, and $\|u\|_{L^\infty(B_1)}$.
\end{lem}
\begin{proof}
We apply Lemma~\ref{lem.hx0} with $U_{x_0} = B_{C_1\gamma}(x_0)\times(-C_2\gamma,C_2\gamma)$ for some constants to be chosen $C_1 \gg C_2$, and study two cases. 

$\bullet$~~Assume $\sup(v-\psi_{x_0})$ is attained at a point $(x_1',y_1)$ (for $x_1'\in \R^{n-1}$, $y\in \R$) on the lateral face of the cylinder $U_{x_0}$, i.e. with $|x_1' - x_0'| = C_1\gamma$ and $0\leq y_1 \leq C_2 \gamma$. Then we have
\begin{align*}
\psi_{x_0} (x_1',y_1) - \varphi(x_1') & \geq  \left(\kappa - \sup |\phi_{\tau \tau}| \right)|x_1'-x_0'|^2- \kappa (n-1)\frac{\Lambda}{\lambda} y_1^2 \\
& \geq  \left(\kappa - \sup |\phi_{\tau \tau}| \right) C_1^2\gamma^2 - \kappa (n-1)\frac{\Lambda}{\lambda} C_2^2 \gamma^2 \geq C_3 \gamma^2,
\end{align*}
provided that $C_1 \gg C_2$. The positive constant $C_3$ depends only on $\kappa$, $n$, the ellipticity constants, $C_1$, and $C_2$. Thus,
\[
v(x_1', y_1) \geq \psi_{x_0} (x_1', y_1) \geq \varphi(x_1')+C_3\gamma^2. 
\]

Now pick a $x_2' \in B_{C_4\gamma}^*(x_1')$ for some positive constant $C_4$ to be chosen and $(x_2'-x_1')\cdot \nabla_{x'}(v-\varphi)(x_1',y_1)\geq 0$. We are considering here $\varphi$ in the whole $B_1$ by simply putting $\varphi(x', y) = \varphi(x')$. Take $\tau = \left( \frac{x_2'-x_1'}{|x_2'-x_1'|},0\right)$, and use the semiconvexity from Lemma~\ref{lem.usim} together with the fact that $\varphi\in C^{1,1}$ to get
\begin{align*}
(v & - \varphi)(x_2',y_1)  = \\
& =  (v -  \varphi)(x_1',y_1) + (x_2'-x_1')\cdot \nabla_{x'} (v - \varphi)(x_1', y_1) + \iint_{[(x_1',y_1),(x_2',y_1)]} (v - \varphi)_{\tau \tau}\\
&  \geq C_3\gamma^2 - C |x_2'-x_1'|^2 \geq ( C_3 - C C_4)\gamma^2 > 0,
\end{align*}
if $C_4$ is chosen appropriately, small enough depending only on $C_3$, $\|\varphi\|_{C^{1,1}}$ and the semiconvexity constant of Lemma~\ref{lem.usim}. Here, and in the next steps, $\iint_{[a,b]}$ denotes the double integral over the segment between the points $a$ and $b$,
\[
\iint_{[a,b]} w := \int_0^{|b-a|} \left[\int_0^s w\left(a+\frac{b-a}{|b-a|} t\right)dt\right]ds.
\] 

To get a contradiction, now suppose that $x_2'\notin S_\gamma$. In particular, this means $v(x_2', 0) = \varphi(x_2')$, and from \eqref{eq.sigmaus} and the semiconcavity in Lemma~\ref{lem.usim} we get 
\begin{align*}
(v- \varphi)(x_2', y_1) & = (v-\varphi)(x_2', 0) + y_1 \frac{\sigma(x_2')}{2} + \iint_{[(x_2',0),(x_2',y_1)]} v_{x_nx_n} \\
& \leq -y_1\frac{\gamma}{2} + Cy_1^2 \leq y_1 \gamma\left(CC_2- \frac{1}{2}\right) \leq 0
\end{align*}
if $C_2$ is small enough depending only on the semiconcavity constant of Lemma~\ref{lem.usim}. Thus, we have reached a contradiction. 

$\bullet$~~Assume now that $\sup(v-\psi_{x_0})$ is attained at a point $(x_1',y_1)$ in the base of the cylinder $U_{x_0}$, i.e. with $|x_1' - x_0'| \leq C_1\gamma$ and $ y_1 = C_2 \gamma$. Then, from $\kappa  > \sup |\varphi_{\tau \tau } |$, we deduce
\[
v(x_1', y_1) \geq \psi_{x_0} (x_1', y_1) \geq \varphi(x_1') - \kappa (n-1) \frac{\Lambda}{\lambda} C_2^2 \gamma^2.
\]

Now choose $x_2'$ such that $|x_2'-x_1'| < C_2\gamma$ and $(x_2'-x_1')\cdot \nabla_{x'}(v-\varphi)(x_1',y_1)\geq 0$. As before, 
\begin{align*}
(v & - \varphi)(x_2',y_1)  = \\
& =  (v -  \varphi)(x_1',y_1) + (x_2'-x_1')\cdot \nabla_{x'} (v - \varphi)(x_1', y_1) + \iint_{[(x_1',y_1),(x_2',y_1)]} (v - \varphi)_{\tau \tau}\\
&  \geq -\kappa (n-1) \frac{\Lambda}{\lambda} C_2^2 \gamma^2 - C |x_2'-x_1'|^2 \geq -C_2^2\left(\kappa (n-1) \frac{\Lambda}{\lambda} + C \right) \gamma^2.
\end{align*}

Now, if $x_2'\notin S_\gamma$ then $v(x_2', 0) = \varphi(x_2')$,
\[
(v- \varphi)(x_2', y_1)  \leq -C_2\frac{\gamma^2}{2} +  \iint_{[(x_2',0),(x_2',y_1)]} v_{x_nx_n}  
 \leq \left(\frac{1}{2}CC_2^2-C_2\right)\gamma^2.
\]
The contradiction follows if one chooses $C_2$ small enough, depending only on $\kappa$, $n$, $\lambda$, $\Lambda$, and the semiconvexity and semiconcavity constants from Lemma~\ref{lem.usim}. 	
\end{proof}

The following lemma is useful to prove the $C^\alpha$ regularity of $\sigma$, and can be found in \cite[Lemma 3.5]{MS08}. It follows from an appropriate use of the strong maximum principle for $\mathcal{M}^-$, the Pucci's extremal operator.

\begin{lem}[\cite{MS08}]
\label{lem.slm}
Let $w$ be a non-negative continuous function in $B_1^*\times (0,1)$ that solves 
\[
\mathcal{M}^- (D^2 w) \leq 0 ~~\textrm{  in  }~~ B_1^* \times (0,1).
\]
Assume 
\[
\limsup\limits_{x_n\downarrow 0^+} w(x', x_n) \geq 1 ~~\textrm{ for }~~ x' \in B_\delta^*(\bar{x}'),
\]
for some ball $B_\delta^*(\bar{x}')\subset B_1^*$. Then 
\[
w(x) \geq \varepsilon > 0 ~~\textrm{ for }~~ x \in B_{1/2}^* \times \left[\frac{1}{4},\frac{3}{4}\right],
\]
for some $\varepsilon$ depending only on $\delta$, and the ellipticity constants $\lambda$ and $\Lambda$.
\end{lem}

We now show the following lemma, analogous to \cite[Lemma 4]{Caf79} (or \cite[Lemma 3.6]{MS08}).

\begin{lem}
\label{lem.sigma_alpha0}
Let $\sigma$ as defined in \eqref{eq.sigmadef}-\eqref{eq.sigmaus}, for $u$ the solution to the thin obstacle problem \eqref{eq.thinobst}. Let $x_0'\in \Omega^*$, then
\[
\sigma(x') \geq -C\left(\|u\|_{L^\infty(B_1)}+\|\varphi\|_{C^{1,1}(B_1^*)} \right)|x'-x_0'|^\alpha,~{\rm for }~x' \in B_1^*
\] 
for some $\alpha > 0$ and $C$ depending only on $n$, $\lambda$, and $\Lambda$. 
\end{lem}
\begin{proof}
Define
\[
K_0 := \|u\|_{L^\infty(B_1)}+\|\varphi\|_{C^{1,1}(B_1^*)},
\]
and notice that by taking $u/K_0$ instead of $u$ if necessary we can assume
\[
\|u\|_{L^\infty(B_1)}+\|\varphi\|_{C^{1,1}(B_1^*)} \leq 1.
\]
Indeed, if $K_0 \geq 1$ then
\[
F_{K_0}(D^2 u) := \frac{1}{K_0} F(D^2(K_0u)),
\]
is a convex elliptic operator with ellipticity constants $\lambda$ and $\Lambda$, and $u/K_0$ is a solution to the nonlinear thin obstacle problem for the operator $F_{K_0}$ with obstacle $\varphi/{K_0}$. In this case, 
\[
\left\| u/K_0\right\|_{L^\infty(B_1)} +\left\| \varphi/K_0 \right\|_{C^{1,1}(B_1^*)} = 1,
\]
as we wanted to see. Thus, from now on we assume $K_0\leq 1$. 

Using Lemmas \ref{lem.sigmaneg}, \ref{lem.usim}, \ref{lem.sgamma} and \ref{lem.slm}, now the proof of this lemma is very similar to the proof of \cite[Lemma 3.6]{MS08}. We give it here for completeness. 

We will show
\begin{equation}
\label{eq.estimatexi}
\sigma(x') \geq -C|x'-x_0'|^\alpha,
\end{equation}
with $C$ and $\alpha > 0$ depending only on $n$, $\lambda$, and $\Lambda$. 

Recall that $\sigma(x') = 2\lim_{x_n \downarrow 0^+} v_{x_n}(x',x_n)$, and that from Lemma~\ref{lem.usim}, $v_{x_n}$ is bounded and $v_{x_nx_n} \leq C$. Moreover, $\sigma$ is non-positive by Lemma~\ref{lem.sigmaneg}, so that $v_{x_n} \leq Cx_n$ for $x_n > 0$. 

In order to reach \eqref{eq.estimatexi} we will prove $v_{x_n}(x) \geq -\theta^k$ for $x\in B_{\gamma^k}^*(x_0')\times (0,\gamma^k)$. Assume this has been already proved for some $k$ with $0 < \gamma \ll \theta < 1$, and consider the function
\[
w:= \frac{v_{x_n} + \theta^k}{\theta^k-C\mu\gamma^k} \textrm{ in } B_{\mu\gamma^k}^*(x_0')\times (0,\mu\gamma^k)
\]
for $\mu$ small enough. Notice that $w$ fulfils the hypotheses of Lemma~\ref{lem.slm}, so that using it together with Lemma~\ref{lem.sgamma} we get
\[
v_{x_n}(x) \geq -\theta^k + \varepsilon(\theta^k-C\mu \gamma^k)\geq -\theta^k+\frac{1}{2}\varepsilon \theta^k
\]
for $x\in B_{\mu\gamma^k/2}^*(x_0')\times(\mu\gamma^k/4,3\mu\gamma^k/4)$, since $\gamma\ll\theta$.
Now, by means of Lemma~\ref{lem.usim}, $v_{x_nx_n} \leq C$, and therefore, for any $y =(y',y_n) \in B^*_{\mu\gamma^k/2}(x_0')\times(0,\mu\gamma^k/4]$,
\begin{align*}
v_{x_n}(y) & \geq - \int_{y_n}^{\mu\gamma^k/4} v_{x_nx_n}(y', s)ds  +v_{x_nx_n}(y', \mu\gamma^k/4)\\
& \geq -C\left(\frac{\mu\gamma^k}{4} - y_n\right) -\theta^k+\frac{1}{2}\varepsilon \theta^k,
\end{align*}
so that we obtain 
\[
v_{x_n}(x)\geq -\theta^k+ \frac{1}{2}\varepsilon \theta^k - \frac{1}{4}\mu C\gamma^k
\]
for $x\in B^*_{\mu\gamma^k/2}(x_0')\times(0,3\mu\gamma^k/4)$. To end the inductive argument we must see 
\[
\theta^{k+1}\geq \theta^k-\frac{1}{2}\varepsilon\theta^k + \frac{1}{4}\mu C\gamma^k.
\]
For this, we pick $\gamma\ll\theta$ so that the right-hand side is smaller than $(1-\frac{1}{4}\varepsilon)\theta^k$, with $\theta$ larger than $1-\frac{1}{4}\varepsilon$. Then, the inductive argument is completed, and \eqref{eq.estimatexi} follows.
\end{proof}

\subsection{Proof of Theorem~\ref{thm.main.1}}
Before proving our main result, let us show the following compactness lemma.
\begin{lem}
\label{lem.comp}
Let $F$ be a nonlinear operator satisfying \eqref{eq.F}, and let $w$ be a continuous function defined on $B_1$. Suppose that $w$ satisfies the problem
\begin{equation}
  F(D^2w)=0 \textrm{ in } B_1^+\cup B_1^-,
\end{equation}
and that
\[
\|w\|_{L^\infty (B_{1})} = 1,~~~~[w]_{{\rm Lip}(B_{1})} \leq 1.
\]

Let $\psi$ be the solution to 
\begin{equation}
  \left\{ \begin{array}{rcll}
  F(D^2\psi)&=&0& \textrm{ in } B_1\\
  \psi&=&w& \textrm{ on } \de B_1,\\
  \end{array}\right.~~
\end{equation}
and let us define the following operator
\[
\tilde{\sigma}(w) := \lim_{h_n\downarrow 0} \left((\de_{x_n} w)(x',h_n) - (\de_{x_n} w)(x',-h_n)\right).
\]

Then, for every $\varepsilon > 0$ there exists some $\eta = \eta(\varepsilon, n, \lambda, \Lambda)> 0$ such that if 
\[
\|\tilde \sigma(w)\|_{L^\infty(B_1^*)} < \eta
\]
then
\[
\|\psi - w\|_{L^\infty(B_{1})} < \varepsilon,
\]
i.e., $\psi$ approximates $w$ as $\eta$ goes to 0. 
\end{lem}
\begin{proof}
Let us argue by contradiction. Suppose that there exists some fixed $\varepsilon > 0$, a sequence of functions $w_k$ and a sequence of convex nonlinear operators uniformly elliptic with ellipticity constants $\lambda$ and $\Lambda$, $F_k$, with $F_k (0) = 0$, such that 
\begin{equation}
  F_k(D^2w_k)=0 \textrm{ in } B_1^+\cup B_1^-
\end{equation}
and 
\[
\|w_k\|_{L^\infty (B_{1})} = 1,~~~~[w_k]_{{\rm Lip}(B_{1})} \leq 1,
\]
with 
\begin{equation}
\label{eq.sigmak}
\|\tilde \sigma (w_k)\|_{L^\infty(B_1^*)} < \eta_k
\end{equation}
for some sequence $\eta_k \to 0$, but such that
\begin{equation}
\label{eq.sigmak2}
\|\psi_k - w_k\|_{L^\infty(B_{1})} \geq  \varepsilon,
\end{equation}
for all $k$, where $\psi_k$ is the solution to 
\begin{equation}
  \left\{ \begin{array}{rcll}
  F_k(D^2\psi_k)&=&0& \textrm{ in } B_1\\
  \psi_k&=&w_k& \textrm{ on } \de B_1.\\
  \end{array}\right.~~
\end{equation}

By Arzelà-Ascoli, up to a subsequence, $w_k$ converges to some function $\bar w$ uniformly in $B_1$, with $\|\bar w\|_{L^\infty(B_1)} = 1$. On the other hand, since $F_k(0) = 0$ and they are uniformly elliptic and convex, they converge up to subsequences, uniformly over compact sets, to some convex nonlinear operator $\bar F$ uniformly elliptic with ellipticity constants $\lambda$ and $\Lambda$ such that $\bar F(0) = 0$. Notice also that $\psi_k$ converges uniformly to the solution $\bar \psi$ to
\begin{equation}
  \left\{ \begin{array}{rcll}
  \bar F(D^2\bar \psi)&=&0& \textrm{ in } B_1\\
  \bar \psi&=&\bar w& \textrm{ on } \de B_1.\\
  \end{array}\right.~~
\end{equation}
and in the limit we obtain, from \eqref{eq.sigmak2},
\begin{equation}
\label{eq.sigmak3}
\|\bar \psi - \bar w\|_{L^\infty(B_{1})} \geq  \varepsilon > 0.
\end{equation}

Now consider the function $w_k + \eta_k |x_n|$ on $B_1$. From \eqref{eq.sigmak}, $w_k + \eta_k |x_n|$ now has a wedge pointing down in the set  $B_1\cup \{x_n = 0\}$, i.e.,
\[
 \tilde\sigma(w_k + \eta_k |x_n|) \geq \eta_k > 0, ~~\textrm{ in } ~~B_1^*.
\]

Therefore, since $F_k(D^2w_k)=0$ in $B_1^+ \cup B_1^-$, we have that, in the viscosity sense,
\[
F_k(D^2 (w_k + \eta_k |x_n|)) \geq 0, ~~\textrm{ in } ~~B_1.
\]

Now, passing to the limit, noticing that $w_k + \eta_k |x_n|$ converges uniformly to $\bar w$ and using \cite[Proposition 2.9]{CC95}, we immediately reach that, in the viscosity sense,
\[
\bar F(D^2 \bar w) \geq 0, ~~ \textrm{ in } ~~B_1.
\]
Repeating the same argument for $w_k - \eta_k |x_n|$ we reach $\bar F(D^2 \bar w) \leq 0$ in $B_1$, to finally obtain
\[
\bar F(D^2 \bar w) = 0, ~~ \textrm{ in } ~~B_1.
\]
This implies $\bar w = \bar \psi$ in $B_1$, which is a contradiction with \eqref{eq.sigmak3}.
\end{proof}

Using the previous results, we now give the proof of Theorem~\ref{thm.main.1}. 

\begin{proof}[Proof of Theorem~\ref{thm.main.1}]
We separate the proof into three steps. In the first step we prove that the solution $u$ is $C^{1, \alpha}$ around points in $\Omega^*$ by means of Lemmas~\ref{lem.sigma_alpha0} and \ref{lem.comp}. In the second step, we use the result from the first step to deduce that $\sigma$ is $C^\alpha$ in $B_{2/3}^*$, to finally complete the proof in the third step.

As in the proof of Lemma~\ref{lem.sigma_alpha0} we assume
\[
\|u\|_{L^\infty(B_1)} + \|\varphi\|_{C^{1,1}(B_1^*)} \leq 1,
\]
to avoid having this constant on each estimate throughout the proof.
\\[3mm]
{\bf Step 1:} Let us suppose that the origin is a free boundary point. Under these circumstances we will prove that there exist some affine function $L = a+b\cdot x$ such that
\begin{equation}
\label{eq.bound0}
\|u-L\|_{L^\infty(B_r)} \leq Cr^{1+\alpha},~~\textrm{ for all } r \geq 0,
\end{equation}
for some constants $C$ and $\alpha > 0$ depending only on $n$, $\lambda$, and $\Lambda$. To do so, we proceed in the spirit of the proof of \cite[Theorem 2]{Caf89}.

Notice that from Lemma~\ref{lem.sigma_alpha0} we know that there exists $\eta > 0$ such that
\begin{equation}
\label{eq.sigmaeps}
|\sigma(x')| \leq \eta |x'|^\alpha, ~~ \textrm{ for all } x'\in B_1^*. 
\end{equation}
Up to replacing from the beginning $u(x)$ by $u(r_0 x)$ with $r_0 \ll 1$, we can make $\eta$ as small as necessary. The choice of the value of $r_0$, and consequently the magnitude in which the constant $\eta$ is made small, will depend only on $n$, $\lambda$, and $\Lambda$. 

Let us show now that there exists $\rho = \rho (\alpha, n, \lambda, \Lambda) < 1$ and a sequence of affine functions 
\begin{equation}
\label{eq.uLk0}
L_k (x) = a_k + b_k\cdot x
\end{equation}
such that 
\begin{equation}
\label{eq.uLk}
\|u-L_k\|_{L^\infty(B_{\rho^k})} \leq \rho^{k(1+\alpha)},
\end{equation}
and 
\begin{equation}
\label{eq.uLk2}
|a_k - a_{k-1}| \leq C\rho^{k(1+\alpha)},~~~~|b_k - b_{k-1}| \leq C\rho^{k\alpha}
\end{equation}
for some constant $C$ depending only on $n$, $\lambda$, and $\Lambda$.

We proceed by induction, taking $L_0 = 0$. Suppose that the $k$-th step is true, and consider 
\[
w_k(x) = \frac{(u-L_k)(\rho^k x)}{\rho^{k(1+\alpha)}}, \textrm{ for } x \in B_1.
\]

Begin by noticing that 
\[
F_k(D^2 w_k) = 0 \textrm{ in } B_1^+\cup B_1^-
\]
for some operator $F_k$ of the form \eqref{eq.F}. On the other hand, from the induction hypothesis, 
\[
\|w_k\|_{B_1} \leq 1. 
\]
Moreover, if we define 
\[
\sigma_k(x') = \lim_{h\downarrow 0 } \left(\de_{x_n} w_k(x',h) - \de_{x_n} w_k (x',-h)\right), ~~ \textrm{ for } x'\in B_1^*,
\]
then one can check that, from \eqref{eq.sigmaeps},
\[
|\sigma_k(x')| \leq \eta |x'|^\alpha.
\]

We apply now Lemma~\ref{lem.comp}. That is, given $\varepsilon > 0$ small, we can choose $\eta$ small enough such that
\[
\|v_k - w_k\|_{L^\infty(B_{1})} \leq \varepsilon,
\]
where $v_k$ is the solution to 
\begin{equation}
  \left\{ \begin{array}{rcll}
  F_k(D^2 v_k)&=&0& \textrm{ in } B_1\\
  v_k&=&w_k& \textrm{ on } \de B_1.\\
  \end{array}\right.~~
\end{equation}

Notice that, by interior estimates, $v_k$ is $C^{2,\alpha}$ in $B_{1/2}$ with estimates depending only on $n$, $\lambda$, and $\Lambda$. Then, let $l_k$ be the linearisation of $v_k$ around 0, so that up to choosing $\rho$,
\begin{align*}
\|w_k - l_k\|_{L^\infty(B_\rho)} & \leq \|w_k - v_k\|_{L^\infty(B_\rho)}+\|v_k - l_k\|_{L^\infty(B_\rho)}\\
& \leq \varepsilon + C\rho^2 \leq \rho^{1+\alpha},
\end{align*}
where $C$ depends only on $n$, $\lambda$, and $\Lambda$, $\rho$ is chosen small enough depending only on $\alpha$, $n$, $\lambda$, and $\Lambda$ so that $C\rho^2 \leq \frac{1}{2} \rho^{1+\alpha}$, and $\eta$ is chosen so that $\varepsilon \leq \frac{1}{2}\rho^{1+\alpha}$. It is important to remark that the choice of $\eta$ depends only on $n$, $\lambda$, and $\Lambda$.

Now, recalling the definition of $w_k$, we reach
\[
\left\|u - L_k - \rho^{k(1+\alpha)}l_k\left(\frac{\cdot}{\rho^k}\right)\right\|_{L^\infty(B_{\rho^{k+1}})} \leq \rho^{(k+1)(1+\alpha)},
\]
so that the inductive step is concluded by taking
\[
L_{k+1}(x) = L_k(x) + \rho^{k(1+\alpha)}l_k\left(\frac{x}{\rho^k}\right).
\]

By noticing that there are bounds on the coefficients of the linearisation of $v_k$ depending only on $n$, $\lambda$, and $\Lambda$, the inequalities in \eqref{eq.uLk2} are obtained. 

Once one has \eqref{eq.uLk0}, \eqref{eq.uLk}, and \eqref{eq.uLk2}, define $L$ as the limit of $L_k$ as $k \to \infty$ (which exists, by \eqref{eq.uLk2}), and notice that, given any $0 < r = \rho^k$ for some $k \in \N$, then 
\begin{align*}
\|u-L\|_{L^\infty(B_r)} \leq \|u-L_k\|_{L^\infty(B_r)} + \sum_{j \geq k} \|L_{j+1} - L_{j}\|_{L^\infty(B_r)} \leq Cr^{1+\alpha}
\end{align*}
for some $C$ depending only on $n$, $\lambda$, and $\Lambda$; as we wanted. 
\\[3mm]
{\bf Step 2:} In this step we prove that the function $\sigma$ defined in \eqref{eq.sigmadef}-\eqref{eq.sigmaus} is $C^\alpha(B_{2/3}^*)$ for some $\alpha = \alpha(n,\lambda,\Lambda)>0$, and
\begin{equation}
\label{eq.sigmaestimate}
\|\sigma\|_{C^{\alpha}(B_{2/3}^*)} \leq C,
\end{equation}
for some constant $C$ depending only on $n$, $\lambda$, and $\Lambda$.

We already know $\sigma$ is regular in the interior of $\Delta^*$ (by boundary estimates) and $\Omega^*$; respectively the coincidence set and its complement in $B_1^*$. In particular, from the interior estimates $\sigma \equiv 0$ in $\Omega^*$. From Lemma~\ref{lem.sigma_alpha0} we also obtain $C^\alpha$ regularity at points in $\de \Delta^*$. Namely, we have that given $(x_0',0)= x_0 \in \de \Delta^*$, 
\begin{equation}
\label{eq.lemsig}
|\sigma(x')|\leq C |x'-x_0'|^\alpha,\textrm{ for } x'\in B_1^*,
\end{equation}
for some constant $C$ depending only on $n$, $\lambda$, and $\Lambda$. 

Therefore, we only need to check that given $x, y \in \Delta$, $x = (x',0)$, $y = (y',0)$, then there exists some $C$ depending only on $n$, $\lambda$, and $\Lambda$ such that, if $|x-y| = r$,
\[
|\sigma(x')-\sigma(y')| \leq Cr^\alpha.
\]

Let $R := \textrm{dist}(x,\Omega)$ and suppose that $\textrm{dist}(x,\Omega) \leq \textrm{dist}(y,\Omega)$. Let $z = (z',0)$, $z'\in \de\Delta^*$, be such that $\textrm{dist}(x, z) = \textrm{dist}(x, \Omega)$, and assume that $\lim_{x_n \downarrow 0^+} \nabla u(z',x_n) = 0$ and $\nabla_{x'} \varphi(z') = 0$ by subtracting an affine function if necessary. Notice that we can do so because we already know from the first step that $u$ has a $C^{1, \alpha}$ estimate around $z'$. Let us then separate two cases:

$\bullet$ If $R < 4r$, then using \eqref{eq.sigmaestimate}
\begin{align*}
|\sigma(x')-\sigma(y')|&  \leq |\sigma(x')-\sigma(z')| + |\sigma(y')-\sigma(z')| \\
& \leq C \left(R^\alpha + (R+r)^\alpha\right) \\
& \leq Cr^\alpha.
\end{align*}

$\bullet$ In the case $R \geq 4r$ we need to use known boundary estimates for this fully nonlinear problem and the previous step of the proof. Notice that $x', y'\in B_{R/2}^*(x')\subset B_{R}^*(x')\subset\Delta^*$, and $u$ restricted to $B_{R}^*(x')$ is thus a $C^{1,1}$ function, since $u = \varphi$ there. In particular, we use that under these hypotheses
\[
R^{1+\alpha}[u]_{C^{1,\alpha}(\overline{B_{R/2}^+(x)})} \leq C\left(\textrm{osc}_{B_R^+(x)} u + R^2[\varphi]_{C^{1,1}(B_R^*(x'))}\right);
\]
see, for example, \cite[Proposition 2.2]{MS06}. Now, remember that the gradient of $u$ at $z$ is 0, so that from the previous step using the bound \eqref{eq.bound0} around $z$,
\begin{equation}
\label{eq.supinf}
|u(p)-\varphi(z')|\leq C|p-z|^{1+\alpha}\leq CR^{1+\alpha} \textrm{ for } p \in B_R^+(x).
\end{equation}
In particular, ${\rm osc}_{B_R^+(x)} u \leq CR^{1+\alpha}$, and thus, this yields
\[
[u]_{C^{1,\alpha}(\overline{B_{R/2}^+(x)})} \leq C,
\]
from which \eqref{eq.sigmaestimate} is proved.
\\[0.3cm]
{\bf Step 3:} Our conclusion now follows by repeating Step 1 around every point on $B_1^*$. Notice that in the first step we only used that the origin was a free boundary point to be able to apply Lemma~\ref{lem.sigma_alpha0} in \eqref{eq.sigmaeps}.

Now, given any point $z'\in B_{1/2}^*$, we can consider the function $u_z$ given by 
\[
u_z(x) := u(x) - \sigma(z')(x_n)^+,
\]
where $(x_n)^+$ denotes the positive part of $x_n$. 

Note that this function fulfils the hypotheses of Step 1, in particular, 
\[
|\sigma_z(x')| := \left| \lim_{h\downarrow 0 } \left(\de_{x_n} u_z(x',h) - \de_{x_n} u_z (x',-h)\right)\right| \leq C|x'-z'|^\alpha, ~~ \textrm{ for } x'\in B_1^*,
\]
for some constant $C$ depending only on $n$, $\lambda$, and $\Lambda$. 

By repeating the exact same procedure as in Step 1, we reach that for every point $z\in B_{1/2}\cup \{x_n = 0\}$, and for every $x \in B_1^+$ there exists some $L_z^+$ affine function such that 
\[
|u(x)-L_z^+| \leq C|x-z|^{1+\alpha},
\]
and the same occurs in $B_1^-$ for a possibly different affine function $L_z^-$.
Therefore, in particular, 
\[
\|u\|_{C^{1,\alpha}(B_{1/2}^*)} \leq C
\]
for some $C$ depending only on $n$, $\lambda$, and $\Lambda$. 

To finish the proof, we could now repeat a procedure like the one done in Step 2, or directly notice that solutions to the nonlinear problem with $C^{1, \alpha}$ boundary data are $C^{1,\alpha}$ up to the boundary (see, for example, \cite[Proposition 2.2]{MS06}). 
\end{proof}

We finally give the:

\begin{proof}[Proof of Corollary \ref{cor.main.1}]
It is an immediate consequence of Theorem~\ref{thm.main.1}. Indeed, consider balls of radius $R_0 := {\rm dist}(K, \de D)$ around points on $K\cap\{x_n = 0\}$ and apply Theorem~\ref{thm.main.1}. To cover the rest of $K$ we use interior estimates, and the result follows by noticing that $\|u\|_{L^\infty(D)}\leq \|g\|_{L^\infty(\de D)}+\|\varphi\|_{L^\infty}$ by the maximum principle.
\end{proof}

\section*{Acknowledgements}
The author would like to thank Alessio Figalli and Xavier Ros-Oton for their guidance and useful discussions on the topics of this paper.

\end{document}